\documentclass[12pt,a4paper]{article}
\usepackage{amsfonts}
\usepackage{amssymb}
\usepackage{amsmath}
\usepackage[perpage,symbol*]{footmisc}
\usepackage{CJK}

\newcommand\songti{\CJKfamily{song}}

\newcommand{\xiaoerhao}{\fontsize{18pt}{\baselineskip}\selectfont}

\newcommand{\xiaowuhao}{\fontsize{9pt}{\baselineskip}\selectfont}

\usepackage{indentfirst}
  \usepackage{geometry}
\begin{document}

\title{\bf\xiaoerhao Two Kinds Hybrid Power Mean Involving Two-Term Exponential Sums and Dedekind Sums}
\author{Yu Jinmin, Han Xue, Wang Tingting\textsuperscript{*}\\
{ \small College of Science, Northwest A$\&$F University, Yangling, Shaanxi, P. R. China}\\
{\small  Corresponding author E-mail: ttwang@nwsuaf.edu.cn}}
\date{}
\maketitle \baselineskip 16pt
\begin{center}
\begin{minipage}{120mm}
{\bf\xiaowuhao\songti Abstract:}{\small\xiaowuhao\songti \  The main purpose of this article is using the analytic mathods and the quadratic residual transformation technique, and properties of Dedekind sums to study the calculating problem of two kinds hybrid power mean involving the two-term exponential sums and Dedekind sums, and give two asymptotic formulas for it. This work is a generalization for existing conclusions.}
 {\bf\xiaowuhao\songti Keywords:}\ \ \xiaowuhao\songti  two-term exponential sums; Dedekind sums; the fourth power mean. 

\end{minipage}
\end{center}
\
\

\section*{1. Introduction}

Let\ $q \geq 3$ be a positive integer. For any integers $m$ and $n$, the two-term exponential sums $C(m,n,k,h;q)$ is defined as follows:
$$
C(m,n,k,h;q) = \sum_{a=0}^{q-1} e\left(\frac{ma^k+na^h}{q}\right),
$$
where\ $e(y)=e^{2\pi i y}$ and $i^2=-1$, $k$ and $h$ are both positive integers.
The study of exponential sums begins with the famous Waring problem. The calculation and upper bound estimation of high-order exponential sums play an important role in analytic number theory, attracted attention from many scholars, see references\ [1-4].

The classical Dedekind sums were introduced by the German mathematician R. Dedekind:

{\bf Definition 1.1}
let\ $q$ be a positive integer and $h$ be an integer prime to $q$. The classical Dedekind sums $S(h,q)$ is defined as follows:
\begin{eqnarray*}
S(h,q) = \sum_{a=0}^{q-1}\left(\left(\frac{a}{q}\right)\right)\left(\left(\frac{ah}{q}\right)\right),
\end{eqnarray*}
where
\begin{eqnarray*}
((x))=\left\{
\begin{array}{ll}\displaystyle x-[x]-\frac{1}{2}, &\textrm{ if\ $x$ is not an integer;} \\ \displaystyle 0,
&\textrm{ if\ $x$ is an integer,}
\end{array}\right.
\end{eqnarray*}
noting that $[x]$ means the integer less than or equal to $x$.

The Dedekind sums describes the properties of $\eta$-function under modular transformations, which is closely related to some important sums and functions in number theory. J. B. Conrey [5] conducted in-depth studies on the mean distribution theorem of $S(d,c)$, and some other interesting results see references\ [7-13].

In \ [14], T. T. Wang and X. W. Pan studied the asymptotic propertities of the mean value of the square two-term exponential sums and Dedekind sums, and obtained some results:
\begin{eqnarray*}
\sum_{m=1}^{p-1}\sum_{n=1}^{p-1}\left|C(m,n,3,1;p)\right|^2\cdot S^{2}(m\overline{n},p) =  \frac{5}{144}p^4+O\left(p^3\cdot {\rm{exp}}\left(\frac{3\ln{p}}{\ln\ln{p}}\right)\right),
\end{eqnarray*}

\begin{eqnarray*}
\sum_{m=1}^{p-1}\sum_{n=1}^{p-1}\left|C(m,n,4,2;p)\right|^2\cdot S^{2}(m\overline{n},p) =  \frac{5}{72}p^4+O\left(p^3\cdot {\rm{exp}}\left(\frac{3\ln{p}}{\ln\ln{p}}\right)\right).
\end{eqnarray*}

In this paper, we study the hybrid mean value of the fourth two-term exponential sums and Dedekind sums with different $k$ and $h$, and give two asymptotic formulas as follows:

{\bf Theorem 1.1.} Let $p>3$ be an odd prime. Then we have the asymptotic formula
\begin{eqnarray*}
&&\sum_{m=1}^{p-1}\sum_{n=1}^{p-1}\left|C(m,n,4,2;p)\right|^4\cdot S^2(m\overline{n},p) = \left\{
\begin{array}{ll}\displaystyle \frac{35}{144}p^5+O\left(p^4\cdot {\rm{exp}}\left(\frac{3\ln{p}}{\ln\ln{p}}\right)\right), &\textrm{ if\ $p\equiv 3\bmod 4$;} \\ \\
\displaystyle \frac{35}{144}p^5+O\left(p^{\frac{9}{2}}\cdot {\rm{exp}}\left(\frac{3\ln{p}}{\ln\ln{p}}\right)\right), &\textrm{ if\ $p\equiv 1\bmod 4$. }
\end{array}\right.
\end{eqnarray*}

{\bf Theorem 1.2.} Let $p>3$ be an odd prime with $p\equiv 3\bmod 4$. Then we have the asymptotic formula
\begin{eqnarray*}
&&\sum_{m=1}^{p-1}\sum_{n=1}^{p-1}\left|C(m,n,5,1;p)\right|^4\cdot S^2(m\overline{n},p) = \frac{5}{48}p^5+O\left(p^4\cdot {\rm{exp}}\left(\frac{3\ln{p}}{\ln\ln{p}}\right)\right).
\end{eqnarray*}

For general integer $q\geq 3$, whether there is a computational formula for
\begin{eqnarray*}
&&\sum_{m=1}^{p-1}\sum_{n=1}^{p-1}\left|C(m,n,k,h;q)\right|^4\cdot S^{2r}(m\overline{n},q)
\end{eqnarray*}
is an open problem, where $k\geq 3$ and $h\geq 3$ are two integers.

\section*{ 2. Several Lemmas}

In this section, as a preparation for proving the theorems, we will give some lemmas.
 \allowdisplaybreaks

{\bf Lemma 2.1.} Let\ $p>3$ be an odd prime, $\lambda$ be a fourth-order character mod $p$, we can get the following identity
\begin{eqnarray*}
&&\sum_{m=1}^{p-1}\sum_{n=1}^{p-1}\left|C(m,n,4,2;p)\right|^4 = \left\{
\begin{array}{ll}\displaystyle 7p^4-18p^3+11p^2, &\textrm{ if\ $p\equiv 3\bmod 4$;} \\ \\
\displaystyle 7p^4-34p^3+p^2(27-4\alpha^2)-4p\alpha^2, &\textrm{ if\ $p\equiv 1\bmod 8$; } \\ \\
\displaystyle  7p^4-26p^3-p^2(19-4\alpha^2)+4p\alpha^2, &\textrm{ if\ $p\equiv 5\bmod 8$. }
\end{array}\right.
\end{eqnarray*}
Note that $\alpha=\sum_{a=1}^{(p-1)/2}\left(\frac{a+\overline{a}}{p}\right)$.

{\bf Proof.} Using the definition of exponential sums and properties of completely residue systems, the fourth two-term exponential sums is expanded as
 \begin{eqnarray}
&&\sum_{m=1}^{p-1}\sum_{n=1}^{p-1}\left|C(m,n,4,2;p)\right|^4\nonumber\\
&=&\sum_{a=0}^{p-1}\sum_{b=0}^{p-1}\sum_{c=0}^{p-1}\sum_{d=0}^{p-1} \sum_{m=1}^{p-1}e\left(\frac{m(a^4+b^4-c^4-d^4)}{p}\right)\sum_{n=1}^{p-1}e\left(\frac{n(a^2+b^2-c^2-d^2)}{p}\right)\nonumber\\
&=&\sum_{a=0}^{p-1}\sum_{b=0}^{p-1}\sum_{c=0}^{p-1}\sum_{d=0}^{p-1} \sum_{m=0}^{p-1}e\left(\frac{m(a^4+b^4-c^4-d^4)}{p}\right)\sum_{n=0}^{p-1}e\left(\frac{n(a^2+b^2-c^2-d^2)}{p}\right)+p^4\nonumber\\
&-&\sum_{a=0}^{p-1}\sum_{b=0}^{p-1}\sum_{c=0}^{p-1}\sum_{d=0}^{p-1}\sum_{m=0}^{p-1}e\left(\frac{m(a^4+b^4-c^4-d^4)}{p}\right) -\sum_{a=0}^{p-1}\sum_{b=0}^{p-1}\sum_{c=0}^{p-1}\sum_{d=0}^{p-1}\sum_{n=0}^{p-1}e\left(\frac{n(a^2+b^2-c^2-d^2)}{p}\right)\nonumber\\
&=&p^2\mathop{\mathop{\sum_{a=0}^{p-1}\sum_{b=0}^{p-1}\sum_{c=0}^{p-1}\sum_{d=0}^{p-1}}_{a^4+b^4-c^4-d^4\equiv 0\bmod p}}_{a^2+b^2-c^2-d^2\equiv 0\bmod p}1+p^4-\sum_{a=0}^{p-1}\sum_{b=0}^{p-1}\sum_{c=0}^{p-1}\sum_{d=0}^{p-1}\sum_{m=0}^{p-1}e\left(\frac{m(a^4+b^4-c^4-d^4)}{p}\right)\nonumber\\ &-&\sum_{a=0}^{p-1}\sum_{b=0}^{p-1}\sum_{c=0}^{p-1}\sum_{d=0}^{p-1}\sum_{n=0}^{p-1}e\left(\frac{n(a^2+b^2-c^2-d^2)}{p}\right)\nonumber\\
&=&p^2\mathop{\mathop{\sum_{a=0}^{p-1}\sum_{b=0}^{p-1}\sum_{c=0}^{p-1}\sum_{d=0}^{p-1}}_{a^4-c^4\equiv d^4-b^4\bmod p}}_{a^2-c^2\equiv d^2-b^2\bmod p}1+p^4-\sum_{a=0}^{p-1}\sum_{b=0}^{p-1}\sum_{c=0}^{p-1}\sum_{d=0}^{p-1}\sum_{m=0}^{p-1}e\left(\frac{m(a^4+b^4-c^4-d^4)}{p}\right)\nonumber\\ &-&\sum_{a=0}^{p-1}\sum_{b=0}^{p-1}\sum_{c=0}^{p-1}\sum_{d=0}^{p-1}\sum_{n=0}^{p-1}e\left(\frac{n(a^2+b^2-c^2-d^2)}{p}\right)\nonumber\\
&=&p^2\mathop{\mathop{\sum_{a=0}^{p-1}\sum_{b=0}^{p-1}\sum_{c=0}^{p-1}\sum_{d=0}^{p-1}}_{a^2+c^2\equiv d^2+b^2\bmod p}}_{a^2-c^2\equiv d^2-b^2 \not\equiv 0\bmod p}1+p^2\mathop{\sum_{a=0}^{p-1}\sum_{b=0}^{p-1}\sum_{c=0}^{p-1}\sum_{d=0}^{p-1}}_{a^2-c^2\equiv d^2-b^2\equiv0 \bmod p}1
-\sum_{a=0}^{p-1}\sum_{b=0}^{p-1}\sum_{c=0}^{p-1}\sum_{d=0}^{p-1}\sum_{m=0}^{p-1}e\left(\frac{m(a^4+b^4-c^4-d^4)}{p}\right)\nonumber\\ &-&\sum_{a=0}^{p-1}\sum_{b=0}^{p-1}\sum_{c=0}^{p-1}\sum_{d=0}^{p-1}\sum_{n=0}^{p-1}e\left(\frac{n(a^2+b^2-c^2-d^2)}{p}\right)+ p^4\nonumber\\
&=&p^2\mathop{\mathop{\sum_{a=0}^{p-1}\sum_{b=0}^{p-1}\sum_{c=0}^{p-1}\sum_{d=0}^{p-1}}_{a^2+c^2\equiv d^2+b^2\bmod p}}_{a^2-c^2\equiv d^2-b^2\bmod p}1+p^2\mathop{\sum_{a=0}^{p-1}\sum_{b=0}^{p-1}\sum_{c=0}^{p-1}\sum_{d=0}^{p-1}}_{a^2-c^2\equiv d^2-b^2\equiv0 \bmod p}1-p^2\mathop{\mathop{\sum_{a=0}^{p-1}\sum_{b=0}^{p-1}\sum_{c=0}^{p-1}\sum_{d=0}^{p-1}}_{a^2+c^2\equiv d^2+b^2\bmod p}}_{a^2-c^2\equiv d^2-b^2\equiv 0\bmod p}1+p^4\nonumber\\
&-&\sum_{a=0}^{p-1}\sum_{b=0}^{p-1}\sum_{c=0}^{p-1}\sum_{d=0}^{p-1}\sum_{m=0}^{p-1}e\left(\frac{m(a^4+b^4-c^4-d^4)}{p}\right) -\sum_{a=0}^{p-1}\sum_{b=0}^{p-1}\sum_{c=0}^{p-1}\sum_{d=0}^{p-1}\sum_{n=0}^{p-1}e\left(\frac{n(a^2+b^2-c^2-d^2)}{p}\right)\nonumber\\
\end{eqnarray}
In (1), we let \begin{eqnarray*}
&&W=\mathop{\mathop{\sum_{a=0}^{p-1}\sum_{b=0}^{p-1}\sum_{c=0}^{p-1}\sum_{d=0}^{p-1}}_{a^2+c^2\equiv d^2+b^2\bmod p}}_{a^2-c^2\equiv d^2-b^2\bmod p}1, N=\mathop{\sum_{a=0}^{p-1}\sum_{b=0}^{p-1}\sum_{c=0}^{p-1}\sum_{d=0}^{p-1}}_{a^2-c^2\equiv d^2-b^2\equiv0 \bmod p}1,S=\mathop{\mathop{\sum_{a=0}^{p-1}\sum_{b=0}^{p-1}\sum_{c=0}^{p-1}\sum_{d=0}^{p-1}}_{a^2+c^2\equiv d^2+b^2\bmod p}}_{a^2-c^2\equiv d^2-b^2\equiv 0\bmod p}1,\\
&&T=\sum_{a=0}^{p-1}\sum_{b=0}^{p-1}\sum_{c=0}^{p-1}\sum_{d=0}^{p-1}\sum_{m=0}^{p-1}e\left(\frac{m(a^4+b^4-c^4-d^4)}{p}\right), M=\sum_{a=0}^{p-1}\sum_{b=0}^{p-1}\sum_{c=0}^{p-1}\sum_{d=0}^{p-1}\sum_{n=0}^{p-1}e\left(\frac{n(a^2+b^2-c^2-d^2)}{p}\right),
\end{eqnarray*}
then \begin{eqnarray}(1)=p^2W+p^2N-p^2S+p^4-T-M.\end{eqnarray}
Next, we calculate the values of $W, N, S, T, M$ in (2) separately.

First, obviously the sum\begin{eqnarray*}
W=\mathop{\mathop{\sum_{a=0}^{p-1}\sum_{b=0}^{p-1}\sum_{c=0}^{p-1}\sum_{d=0}^{p-1}}_{a^2+c^2\equiv d^2+b^2\bmod p}}_{a^2-c^2\equiv d^2-b^2\bmod p}1
\end{eqnarray*}
is equivalent to the number of solutions to the congruence equations system
\begin{eqnarray*}
\left\{
\begin{array}{ll}\displaystyle a^2+c^2\equiv d^2+b^2\bmod p, \\
\displaystyle a^2-c^2\equiv d^2-b^2\bmod p,
\end{array}\right.
\end{eqnarray*}
here a, b, and c go through the residual system module $p$, there is $(2p-1)^2$ solutions, so $W=(2p-1)^2$.

Similarly, \begin{eqnarray*}
N=\mathop{\sum_{a=0}^{p-1}\sum_{b=0}^{p-1}\sum_{c=0}^{p-1}\sum_{d=0}^{p-1}}_{a^2-c^2\equiv d^2-b^2\equiv0 \bmod p}1=(2p-1)^2.
\end{eqnarray*}

Then,
\begin{eqnarray*}
S=\mathop{\mathop{\sum_{a=0}^{p-1}\sum_{b=0}^{p-1}\sum_{c=0}^{p-1}\sum_{d=0}^{p-1}}_{a^2+c^2\equiv d^2+b^2\bmod p}}_{a^2-c^2\equiv d^2-b^2\equiv 0\bmod p}1,
\end{eqnarray*}
we know the number of solutions to this congruence equations system
\begin{eqnarray*}
&&\left\{
\begin{array}{ll}\displaystyle a^2+c^2\equiv d^2+b^2\bmod p, \\
\displaystyle a^2-c^2\equiv d^2-b^2\equiv 0\bmod p,
\end{array}\right.\nonumber\\ \\
&\iff&
\left\{
\begin{array}{ll}\displaystyle a^2\equiv d^2\bmod p, \\
\displaystyle a^2-c^2\equiv d^2-b^2\equiv 0\bmod p,
\end{array}\right.\nonumber\\ \\
&\iff&
\left\{
\begin{array}{ll}\displaystyle a^2\equiv d^2\bmod p, \\
\displaystyle a^2\equiv c^2\bmod p, \\
\displaystyle c^2\equiv b^2\bmod p,
\end{array}\right.\nonumber\\ \\
&\iff&a^2\equiv b^2\equiv c^2\equiv d^2\bmod p,
\end{eqnarray*}
and a, b, and c traverse the residual system module $p$. So $S=8p-7$.

Next,
\begin{eqnarray*}
&&T=\sum_{a=0}^{p-1}\sum_{b=0}^{p-1}\sum_{c=0}^{p-1}\sum_{d=0}^{p-1}\sum_{m=0}^{p-1}e\left(\frac{m(a^4+b^4-c^4-d^4)}{p}\right)\nonumber\\
&=&\sum_{m=0}^{p-1}\left(\sum_{a=0}^{p-1}e\left(\frac{ma^4}{p}\right)\right)^2\left(\sum_{b=0}^{p-1}e\left(\frac{-mb^4}{p}\right)\right)^2,
\end{eqnarray*}
and
\begin{eqnarray*}
&&M=\sum_{a=0}^{p-1}\sum_{b=0}^{p-1}\sum_{c=0}^{p-1}\sum_{d=0}^{p-1}\sum_{n=0}^{p-1}e\left(\frac{n(a^2+b^2-c^2-d^2)}{p}\right)\nonumber\\
&=&\sum_{n=0}^{p-1}\left(\sum_{a=0}^{p-1}e\left(\frac{na^2}{p}\right)\right)^2\left(\sum_{b=0}^{p-1}e\left(\frac{-nb^2}{p}\right)\right)^2.
\end{eqnarray*}
Therefore the formula (1) is equivalent to
\begin{eqnarray*}
&&p^2(2p-1)^2+p^2(2p-1)^2-p^2(8p-7)+p^4-T-M\nonumber\\
&=&9p^4-16p^3+9p^2-T-M.
\end{eqnarray*}

To make further calculate about the values of T and M, we classify prime $p$, if $p\equiv 3\bmod 4$, $T=M=p^4+p^3-p^2$, the formula (1) is equivalent to
\begin{eqnarray}
&&9p^4-16p^3+9p^2-T-M\nonumber\\
&=&7p^4-18p^3+11p^2.
\end{eqnarray}

If $p\equiv 1\bmod 8$, then $\lambda(-1)=1$, $\left(\frac{-1}{p}\right)=1$, by using the characters of Gauss sums, we have
\begin{eqnarray*}
&&T=\sum_{m=0}^{p-1}\left(\sum_{a=0}^{p-1}e\left(\frac{ma^4}{p}\right)\right)^2\left(\sum_{b=0}^{p-1}e\left(\frac{-mb^4}{p}\right)\right)^2\nonumber\\
&=&\sum_{m=1}^{p-1}\left(1+\sum_{a=1}^{p-1}(1+\lambda(a)+\left(\frac{a}{p}\right)+\overline{\lambda}(a))e\left(\frac{ma}{p}\right)\right)^2 \left(1+\sum_{b=1}^{p-1}(1+\lambda(b)+\left(\frac{b}{p}\right)+\overline{\lambda}(b))e\left(\frac{-mb}{p}\right)\right)^2\nonumber\\
&+&p^4\nonumber\\
&=&p^4+\sum_{m=1}^{p-1}\left(\overline{\lambda}(m)\tau(\lambda)+\left(\frac{m}{p}\right)\cdot\sqrt{p}+\lambda(m)\tau(\overline{\lambda})\right)^4\nonumber\\
&=&p^4+\sum_{m=1}^{p-1}\left(\left(\frac{m}{p}\right)\tau^2(\lambda)+p+\left(\frac{m}{p}\right)\tau^2(\overline{\lambda})+2\tau(\lambda)\tau(\overline{\lambda}) +2\sqrt{p}\lambda(m)\tau(\lambda) +2\sqrt{p}\overline{\lambda}(m)\tau(\overline{\lambda})\right)^2\nonumber\\
&=&p^4+\sum_{m=1}^{p-1}\left(2\sqrt{p}\alpha\cdot\left(\frac{m}{p}\right)+3p -2\left(\frac{m}{p}\right)\cdot\sqrt{p}+2\sqrt{p}(\lambda(m)\tau(\lambda)+\overline{\lambda}(m)\tau(\overline{\lambda}))\right)^2\nonumber\\
&=&p^4+17p^3+p^2(4\alpha^2-17)-4p\alpha^2.
\end{eqnarray*}
And then, we can deduce
\begin{eqnarray*}
&&M=\sum_{n=0}^{p-1}\left(\sum_{a=0}^{p-1}e\left(\frac{na^2}{p}\right)\right)^2\left(\sum_{b=0}^{p-1}e\left(\frac{-nb^2}{p}\right)\right)^2\nonumber\\
&=&p^4+\sum_{n=1}^{p-1}\left(1+\sum_{a=1}^{p-1}\left(1+\left(\frac{a}{p}\right)\right)e\left(\frac{na}{p}\right)\right)^2 \left(1+\sum_{b=1}^{p-1}\left(1+\left(\frac{b}{p}\right)\right)e\left(\frac{-nb}{p}\right)\right)^2\nonumber\\
&=&p^4+\sum_{n=1}^{p-1}\left(\left(\frac{n}{p}\right)\cdot\sqrt{p}\right)^4\nonumber\\
&=&p^4+p^3-p^2.
\end{eqnarray*}

So if $p\equiv 1\bmod 8$, the formula (1) is equivalent to
\begin{eqnarray}
&&9p^4-16p^3+9p^2-T-M\nonumber\\
&=&7p^4-34p^3+p^2(27-4\alpha^2)-4p\alpha^2.
\end{eqnarray}

If $p\equiv 5\bmod 8$, then $\lambda(-1)=-1$, $\left(\frac{-1}{p}\right)=1$ and $\tau(\lambda)\cdot\tau(\overline{\lambda})=-p$, we have
\begin{eqnarray*}
&&T=\sum_{m=0}^{p-1}\left(\sum_{a=0}^{p-1}e\left(\frac{ma^4}{p}\right)\right)^2\left(\sum_{b=0}^{p-1}e\left(\frac{-mb^4}{p}\right)\right)^2\nonumber\\
&=&\sum_{m=1}^{p-1}\left(1+\sum_{a=1}^{p-1}(1+\lambda(a)+\left(\frac{a}{p}\right)+\overline{\lambda}(a))e\left(\frac{ma}{p}\right)\right)^2 \left(1+\sum_{b=1}^{p-1}(1+\lambda(b)+\left(\frac{b}{p}\right)+\overline{\lambda}(b))e\left(\frac{-mb}{p}\right)\right)^2\nonumber\\
&+&p^4\nonumber\\
&=&\sum_{m=1}^{p-1}\left(\overline{\lambda}(m)\tau(\lambda)+\left(\frac{m}{p}\right)\cdot\sqrt{p}+\lambda(m)\tau(\overline{\lambda})\right)^2 \left(-\overline{\lambda}(m)\tau(\lambda)+\left(\frac{m}{p}\right)\cdot\sqrt{p}-\lambda(m)\tau(\overline{\lambda})\right)^2\nonumber\\
&+&p^4\nonumber\\
&=&p^4+9p^3-p^2(9-4\alpha^2)-4p\alpha^2.
\end{eqnarray*}
And the other sum
\begin{eqnarray*}
&&M=\sum_{n=0}^{p-1}\left(\sum_{a=0}^{p-1}e\left(\frac{na^2}{p}\right)\right)^2\left(\sum_{b=0}^{p-1}e\left(\frac{-nb^2}{p}\right)\right)^2\nonumber\\
&=&p^4+\sum_{n=1}^{p-1}\left(1+\sum_{a=1}^{p-1}(1+\left(\frac{a}{p}\right))e\left(\frac{na}{p}\right)\right)^2 \left(1+\sum_{b=1}^{p-1}(1+\left(\frac{b}{p}\right))e\left(\frac{-nb}{p}\right)\right)^2\nonumber\\
&=&p^4+\sum_{n=1}^{p-1}\left(\left(\frac{n}{p}\right)\cdot\sqrt{p}\right)^4\nonumber\\
&=&p^4+p^3-p^2.
\end{eqnarray*}

 If $p\equiv 5\bmod 8$, the formula (1) is equivalent to
\begin{eqnarray}
&&7p^4-26p^3-p^2(19-4\alpha^2)+4p\alpha^2.
\end{eqnarray}

Combine (3),(4) and (5), Lemma 1.1 is proved.

{\bf Lemma 2.2.} Let\ $p$ be an odd prime, then we have
\begin{eqnarray*}
&&\sum_{m=1}^{p-1}\sum_{n=1}^{p-1}\left|C(m,n,5,1;p)\right|^4 = \left\{
\begin{array}{ll}\displaystyle 3p^4-p^3(8+2\chi_{2}(-1)+4\chi_{2}(-3))\\+p^2(5+2\chi_{2}(-1)+4\chi_{2}(-3))+2p+1, &\textrm{ if\ $5\nmid p-1$;} \\ \\
\displaystyle 3p^4+O(p^3), &\textrm{ if\ $5\mid p-1$. }
\end{array}\right.
\end{eqnarray*}

{\bf Proof.} From H. Zhang and W. P. Zhang \ [4] we know that
 \allowdisplaybreaks
\begin{eqnarray*}
&&\sum_{m=1}^{p-1}\left|C(m,n,5,1;p)\right|^4 = \left\{
\begin{array}{ll}\displaystyle 3p^3-p^2(8+2\chi_{2}(-1)+4\chi_{2}(-3))-3p, &\textrm{ if\ $5\nmid p-1$;} \\ \\
\displaystyle 3p^3+O(p^2), &\textrm{ if\ $5\mid p-1$. }
\end{array}\right.
\end{eqnarray*}

It is clear that $\sum_{m=1}^{p-1}\sum_{n=1}^{p-1}\left|C(m,n,5,1;p)\right|^4=(p-1)\sum_{m=1}^{p-1}\left|C(m,n,5,1;p)\right|^4$, so we have
\begin{eqnarray*}
&&\sum_{m=1}^{p-1}\sum_{n=1}^{p-1}\left|C(m,n,5,1;p)\right|^4 = \left\{
\begin{array}{ll}\displaystyle 3p^4-p^3(8+2\chi_{2}(-1)+4\chi_{2}(-3))\\+p^2(5+2\chi_{2}(-1)+4\chi_{2}(-3))+2p+1, &\textrm{ if\ $5\nmid p-1$;} \\ \\
\displaystyle 3p^4+O(p^3), &\textrm{ if\ $5\mid p-1$. }
\end{array}\right.
\end{eqnarray*}

{\bf Lemma 2.3.} If $p\equiv 3\bmod 4$, then for any odd character $\chi_{1}$, $\chi_{2}$ and $\chi_{1}\chi_{2}\neq\chi_{0}$, we have the identity
\begin{eqnarray*}
&&\sum_{a=0}^{p-1}\sum_{b=0}^{p-1}\sum_{c=0}^{p-1}\sum_{d=0}^{p-1}\overline{\chi_{1}}\overline{\chi_{2}}(a^4+b^4-c^4-d^4) \chi_{1}\chi_{2}(a^2+b^2-c^2-d^2)=0.
\end{eqnarray*}

{\bf Proof.} From the properties of Dirichlet character sums, noting that $\chi^2_{1}\chi^2_{2}$ is not a principal character mod $p$, we have
\begin{eqnarray*}
&&\sum_{a=0}^{p-1}\sum_{b=0}^{p-1}\sum_{c=0}^{p-1}\sum_{d=0}^{p-1}\overline{\chi_{1}}\overline{\chi_{2}}(a^4+b^4-c^4-d^4) \chi_{1}\chi_{2}(a^2+b^2-c^2-d^2)\nonumber\\
&=&\sum_{a=0}^{p-1}\sum_{b=0}^{p-1}\sum_{c=0}^{p-1}\overline{\chi_{1}}\overline{\chi_{2}}(a^4+b^4-c^4)\chi_{1}\chi_{2}(a^2+b^2-c^2)\nonumber\\
&+&\sum_{d=1}^{p-1}\overline{\chi_{1}\chi_{2}}(d^2)\sum_{a=0}^{p-1}\sum_{b=0}^{p-1}\sum_{c=0}^{p-1}\overline{\chi_{1}}\overline{\chi_{2}}(a^4+b^4-c^4-1) \chi_{1}\chi_{2}(a^2+b^2-c^2-1)\nonumber\\
&=&\sum_{a=0}^{p-1}\sum_{c=0}^{p-1}\overline{\chi_{1}}\overline{\chi_{2}}(a^4-c^4)\chi_{1}\chi_{2}(a^2-c^2) +\sum_{b=1}^{p-1}\overline{\chi_{1}\chi_{2}}(b^2)\sum_{a=0}^{p-1}\sum_{c=0}^{p-1}\overline{\chi_{1}}\overline{\chi_{2}}(a^4+1-c^4)\chi_{1}\chi_{2}(a^2+1-c^2)\nonumber\\
&=&\sum_{c=1}^{p-1}\overline{\chi_{1}\chi_{2}}(c^2)\sum_{a=0}^{p-1}\overline{\chi_{1}}\overline{\chi_{2}}(a^4-1)\chi_{1}\chi_{2}(a^2-1)\nonumber\\
&=&0.
\end{eqnarray*}

{\bf Lemma 2.4. } If $p\equiv 3\bmod 4$, then for any odd character $\chi_{1}$, $\chi_{2}$ and  $\chi_{1}\chi_{2}\neq\chi_{0}$, we have the identity
\begin{eqnarray*}
\sum_{a=0}^{p-1}\sum_{b=0}^{p-1}\sum_{c=0}^{p-1}\sum_{d=0}^{p-1}\overline{\chi_{1}}\overline{\chi_{2}}(a^5+b^5-c^5-d^5)\chi_{1}\chi_{2}(a+b-c-d)=0.
\end{eqnarray*}

{\bf Proof. } From the properties of Dirichlet character sums and the method of Lemma 2.3, we can easily conclude
\begin{eqnarray*}
\sum_{a=0}^{p-1}\sum_{b=0}^{p-1}\sum_{c=0}^{p-1}\sum_{d=0}^{p-1}\overline{\chi_{1}}\overline{\chi_{2}}(a^5+b^5-c^5-d^5)\chi_{1}\chi_{2}(a+b-c-d)=0.
\end{eqnarray*}

{\bf Lemma 2.5. } Let $q\geq 2$ be an integer. Then for any integer $a$ with $(a,q)=1$, we have the identity
\begin{eqnarray*}
&&S(a,q)=\frac{1}{\pi^2q}\mathop{\sum}_{d\mid q}\frac{d^2}{\phi(d)}\mathop{\mathop{\sum}_{\chi\bmod d}}_{\chi(-1)=-1}\chi(a)|L(1,\chi)|^2,
\end{eqnarray*}
where $L(1, \chi)$ denotes the Dirichlet L-function corresponding to character mod $d$.

{\bf Proof. } See [6, Lemma 2].

{\bf Lemma 2.6. } Let $q\geq 2$ be an integer, then we have the asymptotic formula
\begin{eqnarray*}
&&\mathop{\mathop{\sum}_{\chi\bmod d}}_{\chi(-1)=-1}|L(1,\chi)|^4 = \frac{5\pi^4}{144}\phi(q)\prod_{p\mid q}\frac{(p^2-1)^3}{p^4(p^2+1)}+O\left(\frac{\phi(q)}{q}\rm{exp}\left(\frac{3\ln q}{\ln\ln q}\right)\right),
\end{eqnarray*}
where exp(y)=$e^y$, $\prod_{p\mid q}$ denotes the product overall distinct prime divisors of $q$.

{\bf Proof. } See [6, Lemma 3].

{\bf Lemma 2.7. } From the properties of Dirichlet character sums we have
\begin{eqnarray*}
&&\sum_{a=0}^{p-1}\sum_{b=0}^{p-1}\sum_{c=0}^{p-1}\sum_{d=0}^{p-1} \left(\frac{a^4+b^4-c^4-d^4}{p}\right)\left(\frac{a^2+b^2-c^2-d^2}{p}\right)=O(p^{\frac{5}{2}}).\nonumber\\
\end{eqnarray*}

{\bf Proof. } According to the properties of the complete residue system mod $p$, we are aware of that if $a$, $b$ and $c$ pass through a complete residue system mod $p$ respectively, then $ab$, $a+c$ and $b+1$ also pass through the complete residue system mod $p$ respectively. After some transformation for $a$, $b$ and $c$, we have Subsequently, we can conclude that
\begin{eqnarray}
&&\sum_{a=0}^{p-1}\sum_{b=0}^{p-1}\sum_{c=0}^{p-1}\sum_{d=0}^{p-1} \left(\frac{a^4+b^4-c^4-d^4}{p}\right)\left(\frac{a^2+b^2-c^2-d^2}{p}\right)\nonumber\\
&=&\sum_{a=0}^{p-1}\sum_{b=0}^{p-1}\sum_{c=0}^{p-1}\sum_{d=1}^{p-1}\left(\frac{d^6}{p}\right)\left(\frac{a^4+b^4-c^4-1}{p}\right)\left(\frac{a^2+b^2-c^2-1}{p}\right)\nonumber\\ &&+\sum_{a=0}^{p-1}\sum_{b=1}^{p-1}\sum_{c=0}^{p-1}\left(\frac{a^4+b^4-c^4}{p}\right)\left(\frac{a^2+b^2-c^2}{p}\right) +\sum_{a=0}^{p-1}\sum_{c=0}^{p-1}\left(\frac{a^4-c^4}{p}\right)\left(\frac{a^2-c^2}{p}\right)\nonumber\\
&=&(p-1)\sum_{a=0}^{p-1}\sum_{b=0}^{p-1}\sum_{c=0}^{p-1}\left(\frac{a^4+b^4-c^4-1}{p}\right)\left(\frac{a^2+b^2-c^2-1}{p}\right)\nonumber\\ &&+(p-1)\sum_{a=0}^{p-1}\sum_{c=0}^{p-1}\left(\frac{a^4-c^4+1}{p}\right)\left(\frac{a^2-c^2+1}{p}\right) +\sum_{a=0}^{p-1}\sum_{c=0}^{p-1}\left(\frac{a^2+c^2}{p}\right)\nonumber\\
&=&(p-1)\sum_{a=0}^{p-1}\sum_{b=0}^{p-1}\sum_{c=0}^{p-1}\left(\frac{a^4+b^4-c^4-1}{p}\right)\left(\frac{a^2+b^2-c^2-1}{p}\right)\nonumber\\ &&+(p-1)\sum_{a=0}^{p-1}\sum_{c=0}^{p-1}\left(\frac{a^4-c^4+1}{p}\right)\left(\frac{a^2-c^2+1}{p}\right) +1.
\end{eqnarray}

In (6), we calculate
\begin{eqnarray}
&&\sum_{a=0}^{p-1}\sum_{b=0}^{p-1}\sum_{c=0}^{p-1}\left(\frac{a^4+b^4-c^4-1}{p}\right) \left(\frac{a^2+b^2-c^2-1}{p}\right)\nonumber\\
&=&\sum_{a=0}^{p-1}\sum_{b=0}^{p-1}\sum_{c=0}^{p-1}\left(\frac{(a+c)^4+(b+1)^4-c^4-1}{p}\right) \left(\frac{(a+c)^2+(b+1)^2-c^2-1}{p}\right)\nonumber\\
&=&\sum_{a=0}^{p-1}\sum_{b=1}^{p-1}\sum_{c=0}^{p-1}\left(\frac{a^4+4a^3c+6a^2c^2+4ac^3+b^4+4b^3+6b^2+4b}{p}\right) \left(\frac{a^2+2ac+b^2+2b}{p}\right)\nonumber\\
&&+\sum_{a=0}^{p-1}\sum_{c=0}^{p-1}\left(\frac{a^4+4a^3c+6a^2c^2+4ac^3}{p}\right) \left(\frac{a^2+2ac}{p}\right)\nonumber\\
&=&\sum_{a=0}^{p-1}\sum_{b=1}^{p-1}\sum_{c=0}^{p-1}\left(\frac{a^4+4a^3c+6a^2c^2+4ac^3+1+4b+6b^2+4b^3}{p}\right) \left(\frac{a^2+2ac+1+2b}{p}\right)\nonumber\\
&&+\sum_{a=0}^{p-1}\sum_{c=1}^{p-1}\left(\frac{a^4+4a^3c+6a^2c^2+4ac^3}{p}\right) \left(\frac{a^2+2ac}{p}\right)+p\nonumber\\
&=&\sum_{a=0}^{p-1}\sum_{b=1}^{p-1}\sum_{c=0}^{p-1}\left(\frac{2}{p}\right) \nonumber\\ &&\left(\frac{2a^4+a(8c^3+12ac^2+6a^2c+a^3+a^2(2c+a)-2a^3)+8b^3+12b^2+6b+1+2b+1}{p}\right)\nonumber\\ &&\left(\frac{a(2c+a)+1+2b}{p}\right)+(p-1)\sum_{a=0}^{p-1}\left(\frac{a^4+4a^3+6a^2+4a}{p}\right) \left(\frac{a^2+2a}{p}\right)+p\nonumber\\
&=&\sum_{a=0}^{p-1}\sum_{b=0}^{p-1}\sum_{c=0}^{p-1}\left(\frac{2}{p}\right) \nonumber\\ &&\left(\frac{2a^4+a(8c^3+12ac^2+6a^2c+a^3+a^2(2c+a)-2a^3)+8b^3+12b^2+6b+1+2b+1}{p}\right)\nonumber\\ &&\left(\frac{a(2c+a)+1+2b}{p}\right)+(p-1)\sum_{a=0}^{p-1}\left(\frac{(a+1)^4-1}{p}\right) \left(\frac{(a+1)^2-1}{p}\right)+p\nonumber\\
&&-\sum_{a=0}^{p-1}\sum_{c=0}^{p-1}\left(\frac{2}{p}\right) \left(\frac{2a^4+a(8c^3+12ac^2+6a^2c+a^3+a^2(2c+a)-2a^3)+2}{p}\right)\left(\frac{a(2c+a)+1}{p}\right)\nonumber\\
&=&\sum_{a=0}^{p-1}\sum_{b=0}^{p-1}\sum_{c=0}^{p-1}\left(\frac{2}{p}\right) \left(\frac{2a^4+a((2c+a)^3+a^2(2c+a)-2a^3)+(2b+1)^3+2b+1}{p}\right)\nonumber\\ &&\left(\frac{a(2c+a)+1+2b}{p}\right)+(p-1)\sum_{a=0}^{p-1}\left(\frac{a^4-1}{p}\right) \left(\frac{a^2-1}{p}\right)+p\nonumber\\
&&-\sum_{a=0}^{p-1}\sum_{c=0}^{p-1} \left(\frac{a^4-c^4+1}{p}\right)\left(\frac{a^2-c^2+1}{p}\right).
\end{eqnarray}
We know that if $a$ $b$ and $c$ go through a complete residue system mod $p$ respectively, $2c+a$ and $2b+1$ also go through the same
complete residue system mod $p$. So in (7), we use $c$ and $b$ instead of $2c+a$ and $2b+1$, (7) equals
\begin{eqnarray}
&&\sum_{a=0}^{p-1}\sum_{b=0}^{p-1}\sum_{c=0}^{p-1}\left(\frac{2}{p}\right) \left(\frac{ac^3+a^3c+b^3+b}{p}\right)\left(\frac{ac+b}{p}\right)+(p-1)\sum_{a=0}^{p-1}\left(\frac{a^4-1}{p}\right) \left(\frac{a^2-1}{p}\right)+p\nonumber\\
&&-\sum_{a=0}^{p-1}\sum_{c=0}^{p-1} \left(\frac{a^4-c^4+1}{p}\right)\left(\frac{a^2-c^2+1}{p}\right)
\end{eqnarray}

Use (8), (6) equals
\begin{eqnarray}
&&(p-1)\sum_{a=0}^{p-1}\sum_{b=0}^{p-1}\sum_{c=0}^{p-1}\left(\frac{2}{p}\right) \left(\frac{ac^3+a^3c+b^3+b}{p}\right)\left(\frac{ac+b}{p}\right)\nonumber\\
&+&(p-1)^2\sum_{a=0}^{p-1}\left(\frac{a^4-1}{p}\right) \left(\frac{a^2-1}{p}\right)+p(p-1)+1\nonumber\\
&=&(p-1)\sum_{a=0}^{p-1}\sum_{b=1}^{p-1}\sum_{c=0}^{p-1}\left(\frac{2}{p}\right) \left(\frac{ac^3+a^3c+b^3+b}{p}\right)\left(\frac{ac+b}{p}\right)+(p-1)^2\sum_{a=0}^{p-1}\left(\frac{a^4-1}{p}\right) \left(\frac{a^2-1}{p}\right)\nonumber\\
&&+(p-1)\sum_{a=0}^{p-1}\sum_{c=0}^{p-1}\left(\frac{2}{p}\right) \left(\frac{ac^3+a^3c}{p}\right)\left(\frac{ac}{p}\right)+p(p-1)+1\nonumber\\
&=&(p-1)\sum_{a=0}^{p-1}\sum_{b=1}^{p-1}\sum_{c=0}^{p-1}\left(\frac{2}{p}\right) \left(\frac{ac^3+b^2(a^3c+1)+1}{p}\right)\left(\frac{ac+1}{p}\right) +(p-1)^2\sum_{a=0}^{p-1}\left(\frac{a^4-1}{p}\right)\left(\frac{a^2-1}{p}\right)\nonumber\\
&&+(p-1)\sum_{a=0}^{p-1}\sum_{c=0}^{p-1}\left(\frac{2}{p}\right) \left(\frac{a^2+c^2}{p}\right)+p(p-1)+1\nonumber\\
&=&(p-1)\sum_{a=0}^{p-1}\sum_{b=1}^{p-1}\sum_{c=0}^{p-1}\left(\frac{2}{p}\right) \left(\frac{ac^3+b^2(a^3c+1)+1}{p}\right)\left(\frac{ac+1}{p}\right) +(p-1)^2\sum_{a=0}^{p-1}\left(\frac{a^4-1}{p}\right)\left(\frac{a^2-1}{p}\right)\nonumber\\
&&+(p-1)^2\sum_{a=0}^{p-1}\left(\frac{2}{p}\right) \left(\frac{a^2+1}{p}\right)+(p-1)^2\left(\frac{2}{p}\right)+p(p-1)+1
\end{eqnarray}

For any integer $n$, note that the identity
\begin{eqnarray*}
&&\sum_{a=0}^{p-1}\left(\frac{a^2+n}{p}\right)=\left\{
\begin{array}{ll}\displaystyle p-1, &\textrm{ if\ $(n,p)=p$;} \\ \\
\displaystyle -1, &\textrm{ if\ $(n,p)=1$, }
\end{array}\right.
\end{eqnarray*}
(9) equals
\begin{eqnarray*}
&&p(p-1)\mathop{\sum_{a=0}^{p-1}\sum_{c=0}^{p-1}}_{ac^3+1\equiv 0\bmod p}\left(\frac{2}{p}\right) \left(\frac{a^3c+1}{p}\right)\left(\frac{ac+1}{p}\right) -(p-1)\sum_{a=0}^{p-1}\sum_{c=0}^{p-1}\left(\frac{2}{p}\right) \left(\frac{a^3c+1}{p}\right)\left(\frac{ac+1}{p}\right) \nonumber\\
&&-(p-1)\sum_{a=0}^{p-1}\sum_{c=0}^{p-1}\left(\frac{2}{p}\right) \left(\frac{ac^3+1}{p}\right)\left(\frac{ac+1}{p}\right) +(p-1)^2\sum_{a=0}^{p-1}\left(\frac{a^4-1}{p}\right)\left(\frac{a^2-1}{p}\right) \nonumber\\ &&-2(p-1)^2\left(\frac{2}{p}\right)+p(p-1)+1\nonumber\\
&=&p(p-1)\mathop{\sum_{a=0}^{p-1}\sum_{c=0}^{p-1}}_{a\equiv -\overline{c}^3\bmod p}\left(\frac{2}{p}\right) \left(\frac{a^3c+1}{p}\right)\left(\frac{ac+1}{p}\right) -2(p-1)\sum_{a=0}^{p-1}\sum_{c=0}^{p-1}\left(\frac{2}{p}\right) \left(\frac{a^2c+1}{p}\right)\left(\frac{c+1}{p}\right) \nonumber\\
&&+(p-1)^2\sum_{a=0}^{p-1}\left(\frac{a^4-1}{p}\right)\left(\frac{a^2-1}{p}\right) -2(p-1)^2\left(\frac{2}{p}\right)+p(p-1)+1\nonumber\\
&=&p(p-1)\sum_{c=0}^{p-1}\left(\frac{2}{p}\right) \left(\frac{c^8-1}{p}\right)\left(\frac{c^2-1}{p}\right) -2(p-1)\sum_{a=0}^{p-1}\sum_{c=0}^{p-1}\left(\frac{2}{p}\right) \left(\frac{a^2(c^2+c)+c+1}{p}\right) \nonumber\\
&&+(p-1)^2\sum_{a=0}^{p-1}\left(\frac{a^4-1}{p}\right)\left(\frac{a^2-1}{p}\right) -2(p-1)^2\left(\frac{2}{p}\right)+p(p-1)+1\nonumber\\
\end{eqnarray*}

From Weil's classic work \ [15], we know that if $\chi$ is a $q$-th order character to the prime modulo $p$,
and if polynomial $f(x)$ is not a perfect $q$-th power modulo $p$, we have the estimate
\begin{eqnarray*}
&&\sum_{x=N+1}^{N+H}\chi(f(x))=O(p^{\frac{1}{2}}\ln p)
\end{eqnarray*}
where $N$ and $H$ are any positive integers.

So from the estimate, we have
\begin{eqnarray*}
&&p(p-1)\sum_{c=0}^{p-1}\left(\frac{2}{p}\right) \left(\frac{c^8-1}{p}\right)\left(\frac{c^2-1}{p}\right)=O(p^{\frac{5}{2}}).
\end{eqnarray*}

Similarly, we can get the formula
\begin{eqnarray*}
&&(p-1)\sum_{a=0}^{p-1}\sum_{c=0}^{p-1}\left(\frac{2}{p}\right) \left(\frac{a^2(c^2+c)+c+1}{p}\right)\nonumber\\
&=&p(p-1)\mathop{\sum_{c=0}^{p-1}}_{c+1\equiv 0\bmod p}\left(\frac{2}{p}\right) \left(\frac{c^2+c}{p}\right)-(p-1)\sum_{c=0}^{p-1}\left(\frac{2}{p}\right) \left(\frac{c^2+c}{p}\right)\nonumber\\
&=&0
\end{eqnarray*}
so we have
\begin{eqnarray*}
&&\sum_{a=0}^{p-1}\sum_{b=0}^{p-1}\sum_{c=0}^{p-1}\sum_{d=0}^{p-1} \left(\frac{a^4+b^4-c^4-d^4}{p}\right)\left(\frac{a^2+b^2-c^2-d^2}{p}\right)=O(p^{\frac{5}{2}}).\nonumber\\
\end{eqnarray*}

This proves Lemma 2.7.

\section*{ 3. Proof of Theorem 1.1}

We prove the Theorem 1.1 by substituting Lemma 2.5 and Lemma 2.6 into $S(a,q)$ and $L(1, \chi)$. We have
\begin{eqnarray*}
S(a,p)=\frac{1}{\pi^2}\cdot\frac{p}{p-1}\cdot\mathop{\mathop{\sum}_{\chi\bmod p}}_{\chi(-1)=-1}\chi(a)|L(1,\chi)|^2.
\end{eqnarray*}

Using the properties of $C(m,n,k,h;q)$ and combining Lemma 2.1 we have
 \allowdisplaybreaks
\begin{eqnarray}
&&\sum_{m=1}^{p-1}\sum_{n=1}^{p-1}|C(m,n,4,2;p)|^4\cdot S^2(m\overline{n},p)\nonumber\\
&=&\frac{p^2\cdot\pi^{-4}}{(p-1)^2}\mathop{\mathop{\sum}_{\chi_{1}\bmod p}}_{\chi_{1}(-1)=-1}\mathop{\mathop{\sum}_{\chi_{2}\bmod p}}_{\chi_{2}(-1)=-1}|L(1,\chi_{1})|^2\cdot|L(1,\chi_{2})|^2\times\sum_{m=1}^{p-1}\sum_{n=1}^{p-1} \chi_{1}\chi_{2}(m\overline{n})|C(m,n,4,2;p)|^4\nonumber\\
&=&\frac{p^2\cdot\pi^{-4}}{(p-1)^2}\mathop{\mathop{\mathop{\sum}_{\chi_{1}\bmod p}}_{\chi_{1}(-1)=-1}\mathop{\mathop{\sum}_{\chi_{2}\bmod p}}_{\chi_{2}(-1)=-1}}_{\chi_{1}\chi_{2}=\chi_{0}} |L(1,\chi_{1})|^2\cdot|L(1,\chi_{2})|^2\times\sum_{m=1}^{p-1}\sum_{n=1}^{p-1}|C(m,n,4,2;p)|^4\nonumber\\
&+&\frac{p^2\cdot\pi^{-4}}{(p-1)^2}\mathop{\mathop{\mathop{\sum}_{\chi_{1}\bmod p}}_{\chi_{1}(-1)=-1}\mathop{\mathop{\sum}_{\chi_{2}\bmod p}}_{\chi_{2}(-1)=-1}}_{\chi_{1}\chi_{2}\neq\chi_{0}} |L(1,\chi_{1})|^2|L(1,\chi_{2})|^2\nonumber\\
&\times&\tau(\chi_{1}\chi_{2})\tau(\overline{\chi_{1}\chi_{2}})\sum_{a=0}^{p-1}\sum_{b=0}^{p-1} \sum_{c=0}^{p-1}\sum_{d=0}^{p-1} \overline{\chi_{1}}\overline{\chi_{2}}(a^4+b^4-c^4-d^4)\chi_{1}\chi_{2}(a^2+b^2-c^2-d^2),
\end{eqnarray}
where $\chi_{0}$ is the principal character modulo $p$.

From Lemma 2.6 we can conclude
\begin{eqnarray*}
&&\mathop{\mathop{\mathop{\sum}_{\chi_{1}\bmod p}}_{\chi(-1)=-1}\mathop{\mathop{\sum}_{\chi_{2}\bmod p}}_{\chi(-1)=-1}}_{\chi_{1}\chi_{2}=\chi_{0}} |L(1,\chi_{1})|^2|L(1,\chi_{2})|^2\nonumber\\
&=&\mathop{\sum_{\chi\bmod p}}_{\chi(-1)=-1}|L(1,\chi)|^4=\frac{5\pi^4}{144}p+O\left({\rm{exp}}\left(\frac{3\ln{p}}{\ln\ln{p}}\right)\right).
\end{eqnarray*}

Therefore, (10) is equivalent to
\begin{eqnarray}
&&\frac{5p^3}{144(p-1)^2}\times\sum_{m=1}^{p-1}\sum_{n=1}^{p-1}|C(m,n,4,2;p)|^4\nonumber\\
&+&\frac{p^2\cdot\pi^{-4}}{(p-1)^2}\mathop{\mathop{\mathop{\sum}_{\chi_{1}\bmod p}}_{\chi_{1}(-1)=-1}\mathop{\mathop{\sum}_{\chi_{2}\bmod p}}_{\chi_{2}(-1)=-1}}_{\chi_{1}\chi_{2}\neq\chi_{0}} |L(1,\chi_{1})|^2|L(1,\chi_{2})|^2\nonumber\\
&\times&\tau(\chi_{1}\chi_{2})\tau(\overline{\chi_{1}\chi_{2}})\sum_{a=0}^{p-1}\sum_{b=0}^{p-1} \overline{\chi_{1}}\overline{\chi_{2}}(a^4+b^4-c^4-d^4)\chi_{1}\chi_{2}(a^2+b^2-c^2-d^2).
\end{eqnarray}

If $p\equiv 3\bmod 4$, combining Lemma 2.1 and Lemma 2.3, we have
\begin{eqnarray*}
&&\sum_{m=1}^{p-1}\sum_{n=1}^{p-1}|C(m,n,4,2;p)|^4\cdot S^2(m\overline{n},p) = \frac{35}{144}p^5+O\left(p^4\cdot {\rm{exp}}\left(\frac{3\ln{p}}{\ln\ln{p}}\right)\right).
\end{eqnarray*}
This proves the theorem 1 if $p\equiv 3\bmod4$.

If $p\equiv 1\bmod 4$, then for any character $\chi_1\chi_2\neq \chi_0$ in (10), note that $\tau(\chi_{1}\chi_{2})\tau(\overline{\chi_{1}\chi_{2}}) = p$ and $\chi^2_{1}\chi^2_{2}=\chi_{0}$ if and only if $\chi_{1}\chi_{2}=\left(\frac{*}{p}\right)$. So in this case, from (10) we have
\begin{eqnarray}
&&\sum_{m=1}^{p-1}\sum_{n=1}^{p-1}|C(m,n,4,2;p)|^4\cdot S^2(m\overline{n},p)\nonumber\\
&=&\frac{p^2\cdot\pi^{-4}}{(p-1)^2}\mathop{\sum_{\chi\bmod p}}_{\chi(-1)=-1}|L(1,\chi)|^4\cdot\sum_{m=1}^{p-1}\sum_{n=1}^{p-1}|C(m,n,4,2;p)|^4\nonumber\\
&+&\frac{p^3\cdot\pi^{-4}}{(p-1)^2}\mathop{\mathop{\mathop{\sum}_{\chi_{1}\bmod p}}_{\chi_{1}(-1)=-1}\mathop{\mathop{\sum}_{\chi_{2}\bmod p}}_{\chi_{2}(-1)=-1}}_{\chi_{1}\chi_{2}=\left(\frac{*}{p}\right)} |L(1,\chi_{1})|^2\cdot|L(1,\chi_{2})|^2\nonumber\\
&\times&\sum_{a=0}^{p-1}\sum_{b=0}^{p-1}\sum_{c=0}^{p-1}\sum_{d=0}^{p-1}\left(\frac{a^4+b^4-c^4-d^4}{p}\right)\left(\frac{a^2+b^2-c^2-d^2}{p}\right)\nonumber\\
&+&\frac{p^3\cdot\pi^{-4}}{(p-1)^2}\mathop{\mathop{\mathop{\sum}_{\chi_{1}\bmod p}}_{\chi_{1}(-1)=-1}\mathop{\mathop{\sum}_{\chi_{2}\bmod p}}_{\chi_{2}(-1)=-1}}_{\chi_{1}\chi_{2}\neq\chi_{0}, \left(\frac{*}{p}\right)} |L(1,\chi_{1})|^2|L(1,\chi_{2})|^2\nonumber\\
&\times&\sum_{a=0}^{p-1}\sum_{b=0}^{p-1}\sum_{c=0}^{p-1}\sum_{d=0}^{p-1}\overline{\chi_{1}}\overline{\chi_{2}}(a^4+b^4-c^4-d^4)\chi_{1}\chi_{2}(a^2+b^2-c^2-d^2).\nonumber\\
&=&W_{1}+W_{2}+W_{3}.
\end{eqnarray}

From Lemma 2.1 we have
\begin{eqnarray}
&W_{1}&=\frac{p^2\cdot\pi^{-4}}{(p-1)^2}\mathop{\sum_{\chi\bmod p}}_{\chi(-1)=-1}|L(1,\chi)|^4\cdot\sum_{m=1}^{p-1}\sum_{n=1}^{p-1}|C(m,n,4,2;p)|^4\nonumber\\
&&=\frac{35}{144}p^5+O\left(p^4\cdot {\rm{exp}}\left(\frac{3\ln{p}}{\ln\ln{p}}\right)\right).
\end{eqnarray}

From Lemma 2.7 we can conclude $W_{2}= O(p^{\frac{9}{2}})$.

From the method of Lemma 2.3 we have
\begin{eqnarray*}
&W_{3}&=\frac{p^3\cdot\pi^{-4}}{(p-1)^2}\mathop{\mathop{\mathop{\sum}_{\chi_{1}\bmod p}}_{\chi_{1}(-1)=-1}\mathop{\mathop{\sum}_{\chi_{2}\bmod p}}_{\chi_{2}(-1)=-1}}_{\chi_{1}\chi_{2}\neq\chi_{0}, \left(\frac{*}{p}\right)} |L(1,\chi_{1})|^2|L(1,\chi_{2})|^2\nonumber\\
&\times&\sum_{a=0}^{p-1}\sum_{b=0}^{p-1}\sum_{c=0}^{p-1}\sum_{d=0}^{p-1} \overline{\chi_{1}}\overline{\chi_{2}}(a^4+b^4-c^4-d^4)\chi_{1}\chi_{2}(a^2+b^2-c^2-d^2)\nonumber\\
&&=0.
\end{eqnarray*}

So we have the identity
\begin{eqnarray*}
&&\sum_{m=1}^{p-1}\sum_{n=1}^{p-1}|C(m,n,4,2;p)|^4\cdot S^2(m\overline{n},p) = \frac{35}{144}p^5+O\left(p^{\frac{9}{2}}\cdot {\rm{exp}}\left(\frac{3\ln{p}}{\ln\ln{p}}\right)\right).
\end{eqnarray*}

Theorem 1.1 is proved.

\section*{ 4. Proof of Theorem 1.2}

From Lemma 2.2 and the properties of Dedekind sums, we can calculate the formula as follows
\begin{eqnarray}
&&\sum_{m=1}^{p-1}\sum_{n=1}^{p-1}\left|C(m,n,5,1;p)\right|^4\cdot S^2(m\overline{n},p)\nonumber\\
&=&\frac{p^2\cdot\pi^{-4}}{(p-1)^2}\mathop{\mathop{\mathop{\sum}_{\chi_{1}\bmod p}}_{\chi_{1}(-1)=-1}\mathop{\mathop{\sum}_{\chi_{2}\bmod p}}_{\chi_{2}(-1)=-1}}_{\chi_{1}\chi_{2}=\chi_{0}} |L(1,\chi_{1})|^2\cdot|L(1,\chi_{2})|^2\times\sum_{m=1}^{p-1}\sum_{n=1}^{p-1}|C(m,n,5,1;p)|^4\nonumber\\
&+&\frac{p^2\cdot\pi^{-4}}{(p-1)^2}\mathop{\mathop{\mathop{\sum}_{\chi_{1}\bmod p}}_{\chi_{1}(-1)=-1}\mathop{\mathop{\sum}_{\chi_{2}\bmod p}}_{\chi_{2}(-1)=-1}}_{\chi_{1}\chi_{2}\neq\chi_{0}} |L(1,\chi_{1})|^2|L(1,\chi_{2})|^2\nonumber\\
&\times&\tau(\chi_{1}\chi_{2})\tau(\overline{\chi_{1}\chi_{2}})\sum_{a=0}^{p-1}\sum_{b=0}^{p-1} \overline{\chi_{1}}\overline{\chi_{2}}(a^5+b^5-c^5-d^5)\chi_{1}\chi_{2}(a+b-c-d).
\end{eqnarray}

If $p\equiv 3\bmod 4$, combining Lemma 2.2, Lemma 2.4 and (10) we have the asymptotic formula
\begin{eqnarray*}
&&\sum_{m=1}^{p-1}\sum_{n=1}^{p-1}|C(m,n,5,1;p)|^4\cdot S^2(m\overline{n},p) = \frac{5}{48}p^5+O\left(p^4\cdot {\rm{exp}}\left(\frac{3\ln{p}}{\ln\ln{p}}\right)\right).
\end{eqnarray*}

This finishes the proof of the theorem 2.2.

{\bf Authors’ contributions:} All authors have equally contributed to this work. All authors read
and approved the final manuscript.

{\bf Funding:} This work is supported by the N. S. B. R. P. (2022JM-013) of Shaanxi Province.

{\bf Competing interests:} The authors declare that there are no conflicts of interest regarding the
publication of this paper.

 \end{document}